\documentclass{iopart}

\usepackage{t1enc}
\usepackage{graphicx}
\usepackage[english]{babel}
\usepackage[T1]{fontenc}
\usepackage[latin1]{inputenc}
\usepackage{url}
\usepackage{dsfont}
\usepackage{graphicx}
\usepackage{graphics}
\usepackage{subfig}
\usepackage{iopams} 
\usepackage{setstack}
\usepackage{amstext}
\usepackage{amsopn}
\usepackage{amsthm}
\usepackage{enumitem}
\usepackage{color}

\definecolor{MatBlue}{rgb}{0,0.4431,0.7373}
\definecolor{MatOrange}{rgb}{0.9255,0.6902,0.1216}
\definecolor{MatRed}{rgb}{0.8471,0.3216,0.0941}
\definecolor{darkblue}{rgb}{0,0,0.6}

\DeclareMathOperator*{\mydef}{\mathrel{\mathop:}=}

\newcommand{\argminsub}[1]{\underset{{ #1 }}{{\rm argmin}}}

\DeclareMathOperator{\R}{\mathbb{R}}
\DeclareMathOperator{\F}{\mathbb{F}}

\newcommand{\sqnorm}[1]{\Vert #1 \Vert_2^2}
\newcommand{\termabb}[2]{\emph{#1} (\emph{#2})}

\newcommand{\bx}{ {\mathbf x}}
\newcommand{\bu}{ {\mathbf u}}
\newcommand{\A}{\mathcal A}
\newcommand{\M}{\mathcal M}
\newcommand{\Pat}{\mathcal P}
\def\kwave/{\textbf{k}-\texttt{Wave}}

\begin{document}

\title{On the Adjoint Operator in Photoacoustic Tomography}
\author{Simon R.~Arridge$^1$, Marta M.~Betcke$^1$, Ben T.~Cox$^2$, Felix Lucka$^1$ and Brad E.~Treeby$^2$}

\address{$^1$ Department of Computer Science, University College London, WC1E 6BT London, UK}
\address{$^2$ Department of Medical Physics and Bioengineering, University College London, WC1E 6BT London, UK}
\ead{f.lucka@ucl.ac.uk}

\begin{abstract} 
Photoacoustic Tomography (PAT) is an emerging biomedical \textit{imaging from coupled physics} technique, in which the image contrast is due to optical absorption, but the information is carried to the surface of the tissue as ultrasound pulses. Many algorithms and formulae for PAT image reconstruction have been proposed for the case when a complete data set is available. In many practical imaging scenarios, however, it is not possible to obtain the full data, or the data may be sub-sampled for faster data acquisition. In such cases, image reconstruction algorithms that can incorporate prior knowledge to ameliorate the loss of data are required. Hence, recently there has been an increased interest in using  \textit{variational image reconstruction}. A crucial ingredient for the application of these techniques is the adjoint of the PAT forward operator, which is described in this article from physical, theoretical and numerical perspectives. First, a simple mathematical derivation of the adjoint of the PAT forward operator in the continuous framework is presented. Then, an efficient numerical implementation of the adjoint using a $k$-space time domain wave propagation model is described and illustrated in the context of variational PAT image reconstruction, on both 2D and 3D examples including inhomogeneous sound speed. The principal advantage of this analytical adjoint over an algebraic adjoint (obtained by taking the direct adjoint of the particular numerical forward scheme used) is that it can be implemented using currently available fast wave propagation solvers.
\end{abstract}

\submitto{\IP}
\maketitle

\section{Introduction} \label{sec:Intro}

\termabb{Photoacoustic Tomography}{PAT} is a biomedical imaging modality that has come to prominence over the past two decades due to its ability to provide images of soft tissue based on optical absorption with high spatial resolution. 
Optical absorption contrast is desirable because imaging at multiple wavelengths can in principle provide spectroscopic (chemical) information on the absorbing molecules (chromophores). However, due to the highly scattering nature of most biological tissues, imaging purely with light - light in, light out - can only be achieved with high spatial resolution in the vicinity of the tissue surface; beyond that the achievable image resolution falls off quickly with depth. PAT overcomes this problem by generating broadband ultrasonic pulses at the regions of optical absorption which then carry the information about the optical contrast to the tissue surface without being scattered significantly. PAT is increasingly widely used in preclinical studies to image small animal anatomy, pathology and physiology \cite{WaPaKuXiStWa03,YaWa14,XiWa14} and PAT's potential as a clinical technique is also being actively explored \cite{KrLaReDeStDo10,ZaVaGa14,TaNt15}. The spectroscopic information PAT provides (ie. differentiating endogenous chromophores or contrast agents based on their absorption spectra) opens up the possibility of molecular and genomic imaging, and even for functional imaging by monitoring changes in blood oxygenation and flow. For further references on the experimental aspects and applications of PAT see the review papers \cite{Wa09,Bea11,NiXi14}.

\subsection{The Physics of Photoacoustic Tomography} \label{subsec:PhysicsPAT}
PAT is based on the principle that acoustic waves can be generated by the absorption of modulated electromagnetic radiation. When the radiation is in the visible to near-infrared spectrum this is known as the \textit{photoacoustic} or \textit{optoacoustic} effect. The term \textit{thermoacoustic} is used when the radiation is in the microwave region. In PAT experiments, it is common to use pulses of a few nanoseconds duration, such as are generated by a Q-switched laser, to illuminate the tissue. The light will be scattered within the tissue and absorbed by any chromophores present. For the photoacoustic effect to take place, the energy must be subsequently thermalised (as opposed to emitted radiatively) and, if it occurs fast enough, this will give rise to a pressure increase localised to the regions of absorption. These regions of raised pressure and temperature will initiate acoustic and thermal waves as the tissue is elastic and thermally conducting. There are two relevant timescales to be considered. First, for the pressure and temperature fluctuations to be decoupled, the thermal relaxation time must be much longer than the characteristic timescale for the acoustic propagation. This condition is known as \textit{thermal confinement}, and it must be satisfied if the acoustic pressure is to be modelled using a wave equation, rather than a set of coupled equations in temperature and pressure. The second condition, known as \textit{stress} or \textit{pressure confinement}, is that the characteristic timescale for the acoustic propagation, which will depend on the sound speed and the size of the absorbing regions, must be much slower than the time taken to deposit the optical energy as heat. This is essentially an isochoric condition requiring there to be no change is density of the tissue during heat deposition. When this condition is satisfied, the acoustic propagation can be modelled as an initial value problem with the \textit{initial (acoustic) pressure distribution} as the initial condition. The second initial condition - that the time rate of the change of the pressure is zero - derives from the isochoric condition through the assumption of zero initial particle velocity.\\
In PAT, measurements of the induced broadband acoustic pulses are recorded at the surface of the tissue by an array of ultrasound sensors. In order to form a PAT image, it is necessary to solve an inverse initial value problem by inferring an initial acoustic pressure distribution from measured acoustic time series. This image reconstruction step is the subject of this paper. A possible subsequent image reconstruction step, not considered here, is to recover images of the optical properties of the tissue. This requires a model of light transport and is sometimes known as \emph{quantitative PAT} \cite{CoLaArBe12}. The solution of the acoustic inversion examined here provides the data for this optical inversion.

\subsection{Challenges of Photoacoustic Tomography} \label{subsec:ChallengesPAT}

Image reconstruction in PAT calls for the robust inversion of the linear system 
\begin{equation}
f = \A p_0 + \varepsilon, \label{eq:FwdEq}
\end{equation}
where $p_0$ represents the initial pressure distribution we wish to recover, $f$ represents the measured acoustic pressure time series, $\varepsilon$ an additional measurement error, and the linear operator $A$ models both the ultrasonic wave propagation in the tissue and the measurement process (we will discuss the details of this modeling in the Section \ref{subsec:MathModeling}). Solving \eref{eq:FwdEq} can be done in several ways:
\begin{enumerate}
\item Assuming homogeneous acoustic properties and the absence of acoustic absorption, the measured time series can be related to the initial pressure distribution via the spherical mean Radon transform. In this case, the methods of integral geometry can be used to derive direct, explicit inversion formulae for certain sensor geometries, such as e.g. spherical arrays. See \cite{KuKu11} and the references therein for an overview of such methods. 
\item In more general scenarios, a direct inverse operator corresponding to a forward model of the wave equation can be used. One commonly used example is given by \termabb{time reversal}{TR} approaches \cite{FiPa04,XuWa04a,XuWa04b,HrKuNg08}: 
The measured time series at the boundary is used in reversed temporal order as a time dependent Dirichlet boundary condition for the backward wave equation and the pressure field at $t = 0$ is used as an estimate of the initial pressure distribution. With the development of fast numerical methods that are capable of computing full 3D wave propagation with high spatial and temporal resolution, TR and its enhanced variants \cite{StUh09,QiStUhZh11} became an attractive choice in applications \cite{LaNoClZhTrCoJoScLyBe12,HuNiScGuScAnWa12,JaLaOgTrCoZhJoPiPhMaLyPePuBe15}.
\item In particular in scenarios with incomplete or sub-sampled data, numerical models of wave propagation can be used within the \emph{variational image reconstruction} framework to find a regularized least-squares solution of \eref{eq:FwdEq} by solving the optimization problem
\begin{equation}
p_{\mathrm{rec}} = \argminsub{p_0} \left\lbrace \case{1}{2} \sqnorm{\A p_0 -f} + \lambda \mathcal{J}(p_0) \right\rbrace. \label{eq:VarReg}
\end{equation}
Here, $\lambda > 0$ is a \emph{regularization parameter} and $\mathcal{J}(p_0)$ is a suitable \emph{regularization functional} that aims to encode a-priori knowledge about the true solution, $p_0$.
\end{enumerate}
While variational image reconstruction in PAT has recently been shown to yield superior results compared to TR solutions \cite{HuWaNiWaAn13,LuBeArCoHuZhBe15,BeGlSc15}, a potential challenge to its application is that solving \eref{eq:VarReg} by \emph{first order optimization} methods (see, e.g., \cite{BuSaSt14} for an overview) requires an efficient numerical scheme for the application of the \emph{adjoint} PAT operator, $\A^*$. One approach is to first discretise the forward operator and then to find its algebraic adjoint, for example by explicit reversal of the computational steps of the forward scheme \cite{HuWaNiWaAn13}. 
The adjoint obtained by this \textit{discretise-then-adjoint} approach however is not the adjoint of the continuous operator but of the particular numerical scheme. Therefore, in this paper we revisit the mathematical modelling of PAT to derive the general analytical form of the PAT adjoint. The explicit analytical form of the adjoint allows us new theoretical insights, in particular, it clarifies the relation between the adjoint and TR. The analytical adjoint is independent of discretization, i.e., both the forward and adjoint equations can be discretised with any suitable discretization scheme and pre-existing high performance codes can be used to apply the forward and adjoint operators.
The remainder of the paper is structured as follows: In Section \ref{sec:Theory}, we revisit the mathematical modelling of the forward problem in PAT and derive an explicit form of the adjoint PAT operator. Section \ref{sec:MethImp} recapitulates how the adjoint operator can be used within a variational scheme. The utility of \textit{adjoint-then-discretise} paradigm is evaluated in Section \ref{sec:Studies}: After a direct comparison to the \textit{discretise-then-adjoint} paradigm in \cite{HuWaNiWaAn13}, variational methods are compared to approaches based on TR for two image reconstruction problems including sound-speed inhomogeneity. In Section \ref{sec:DisCon}, we summarize our work and discuss its relations to other approaches. An efficient implementation of the adjoint operator using a \textit{k-space pseudospectral wave propagation model} is given in \ref{subsec:ImplOp}.

\section{Theory} \label{sec:Theory}

\subsection{Mathematical Modeling} \label{subsec:MathModeling}

\paragraph{Wave propagation}
Acoustic wave propagation through a compressible fluid is usually modelled mathematically by linearising the equations of fluid dynamics derived from conservation laws and an equation of state. Under the condition that the acoustic particle velocity is much slower than the sound speed, an acoustic wave in tissue can then be modelled by the system of first order equations \cite{TrJaReCo12}:
\numparts
\begin{eqnarray}
\fl \quad \frac{\partial}{\partial t} \bu(\bx,t)& = - \frac{1}{\rho_0(\bx)} \nabla p(\bx,t) \qquad \qquad  \qquad \qquad \quad \,  \text{\small (momentum conservation)}\label{eq:MonCon} \\
\fl \quad \frac{\partial}{\partial t} \rho(\bx,t)& = - \rho_0(\bx) \nabla \cdot \bu(\bx,t) - \bu(\bx,t) \cdot \nabla \rho_0(\bx) + s(\bx,t)  \quad \; \text{\small (mass conservation)} \label{eq:MassCon}\\ 
\fl \quad p(\bx,t)& = c_0^2(\bx) (\rho(\bx,t) + \mathbf{d}\cdot \nabla\rho_0(\bx)), \qquad \qquad \quad  \text{\small (pressure-density relation)} \label{eq:PreDenRel}
\end{eqnarray}
\endnumparts
where $p$ and $\rho$ are the acoustic pressure and acoustic density fluctuations, $\bu$ is the acoustic particle velocity (the time derivative of the acoustic displacement $\mathbf{d}$), $s$ is a mass source term, and $\rho_0$ and $c_0$ are the ambient density and sound speed respectively. Both $\rho_0$ and $c_0$ are positive and bounded from above and below. The system (\ref{eq:MonCon}-\ref{eq:PreDenRel}) can be combined into the second-order wave equation for a heterogeneous medium:
\begin{equation}
\fl \qquad \square^2 p(\bx,t) \mydef \left( \frac{1}{c_0^2(\bx)}\frac{\partial^2}{\partial t^2}  - \rho_0(\bx)\nabla\cdot \left(\frac{1}{\rho_0(\bx)} \nabla \right) \right) p(\bx,t) =  \frac{\partial}{\partial t} s(\bx,t), \label{eq:WaveSecondOrder}
\end{equation}
where we introduced the inhomogenous d'Alembert operator $\square^2$.

\paragraph{Initial and boundary conditions}
As described in Section \ref{subsec:PhysicsPAT},
$p_0$ can be introduced as the initial value of the pressure rather than as a source term. Also, assuming the particle velocity $\bu$ at $t = 0$ vanishes, (\ref{eq:MonCon}-\ref{eq:PreDenRel}) suggests that the time derivative of $p$ should be initialized with zero. With these assumptions, the forward problem of PAT in $d$ dimensions is stated as
\numparts 
\begin{eqnarray}
\square^2 p(\bx,t) & = \; 0, \qquad \qquad &\text{in} \quad {\mathbb{R}}^d \times (0, \infty), \label{eq:fpatA}\\
p(\bx,0)& = \; p_0(\bx), &\text{in} \quad {\mathbb{R}}^d, \label{eq:fpatp0} \\
\frac{\partial p}{\partial t} (\bx,0)& = \; 0, &\text{in} \quad {\mathbb{R}}^d. \label{eq:fpatpt0}
\end{eqnarray}
\endnumparts
Because in practice it is possible to time-gate out any reflections from the experimental equipment, the propagation can be modelled as though it occurs in an unbounded domain i.e., no explicit boundary conditions are required.

\paragraph{Measurement}
We deliberately model the measurement of the pressure field $p(\bx,t)$ in two steps: Firstly, each element of the acoustic sensor array has only limited access to the pressure field: It only detects the field over a small, but finite, volume of space, and the measurement continues only for a finite time $T$. To signify this, we introduce the \textit{window} function $\omega(\bx,t) \in C_0^\infty(\Gamma \times (0,T))$, which maps $p(\bx,t)$ to the field accessible by the sensor array: $g(\bx,t) = \omega(\bx,t) p(\bx,t)$. Here, $\Gamma$ is a $d$-dim, open, bounded set with non-zero measure that contains the cumulated volume of all sensor elements. The sound speed on $\Gamma$ is assumed to be constant. In a second step, each sensor element transforms the accessible spatio-temporal part of the pressure field into a sequence of measured values. This process can be modelled by a measurement operator $\mathcal{M}$ which maps $g(\bx,t) \in C_0^\infty(\Gamma \times (0,T))$ to the data $f \in \R^L$, $L \in \mathbb{N}$. Each type of ultrasonic detector elements commonly used for PAT (see \cite{Bea11,LuRa13} for comprehensive reviews) leads to a different particular form of $\mathcal{M}$. For this work, we only assume that have access to
the adjoint operator $\mathcal{M}^*$.

\subsection{Derivation of the Adjoint PAT Operator} \label{subsec:DerivationAdjoint}

In this section, we derive the adjoint operator $\A^*$ to the forward operator
\numparts 
\begin{eqnarray}
\A &: C_0^\infty (\R^d ) &\rightarrow \; \R^L \label{eq:DefA1} \\
\A & \left[p_0\right] (\bx, t) &= \; \mathcal{M} \; \omega(\bx,t) p(\bx,t), \label{eq:DefA2}
\end{eqnarray}
\endnumparts 
where $p(\bx,t)$ is the time dependent solution of the PAT forward problem (\ref{eq:fpatA}-\ref{eq:fpatpt0}). The assumption $p_0 \in C_0^\infty (\R^d )$, is justified by the short heat diffusion always present in the thermalisation process generating the initial pressure (cf. Section \ref{subsec:PhysicsPAT}). As we assume that $\mathcal{M}^*$ is given, $\A^* = \Pat^* \M^*$, where $\Pat^*: C_0^\infty \left(\Gamma \times (0,T) \right) \rightarrow C_0^\infty (\R^d)$ is the adjoint of 
\numparts 
\begin{eqnarray}
\Pat &: C_0^\infty (\R^d ) &\rightarrow \; C_0^\infty\left(\Gamma \times (0,T)\right) \label{eq:DefP1} \\
\Pat & \left[p_0\right] (\bx, t) &= \; \omega(\bx,t) p(\bx,t), \hspace{4em} \text{in} \quad \Gamma \times (0,T), \label{eq:DefP2}
\end{eqnarray}
\endnumparts 
with respect to the generic $L_2$ bilinear form in $C_0^\infty \left(\Gamma \times (0,T) \right)$ and $C_0^\infty(\mathbb{R}^d)$, respectively. By definition, this adjoint $\Pat^*$ has to satisfy the equality 
\begin{equation}
  \underbrace{\int_0^T \int_{\Gamma} g(\bx, t) \Pat \left[p_0\right](\bx,t)  d\bx dt}_{(\diamond)} = \int_{\R^d}  p_0(\bx) \Pat^* \left[g\right] (\bx) d\bx \label{eq:adj}
\end{equation}
for any $p_0 \in C_0^\infty (\R^d)$, $g \in C_0^\infty \left(\Gamma \times (0,T) \right)$. Therefore, we start with the left hand side $(\diamond)$ of \eref{eq:adj} and we seek to reformulate it in terms of an operator acting on the function $g$. The solution of the initial value problem on an unbounded domain (\ref{eq:fpatA}-\ref{eq:fpatpt0}) can be written as
\begin{equation}\label{eq:Gp0}
p(\bx , t) = \frac{\partial}{\partial t}  \int_{\R^d} G(\bx, t | \bx', 0) p_0(\bx') d\bx',
\end{equation}
where $G(\bx,t|\bx',t')$ is the causal Green's function of the homogeneous problem corresponding to \eref{eq:WaveSecondOrder} \cite{Ba89}:
\begin{equation}
\square^2 G(\bx,t|\bx',t') = \delta_{d}(\bx-\bx')\delta(t-t').
\end{equation}
Substituting $\Pat$, \eref{eq:DefP2}, and $p(\bx,t)$, \eref{eq:Gp0} into $(\diamond)$ we obtain
\begin{equation}\label{eq:adjLHS}
(\diamond) = \int_0^T \int_{\Gamma} g(\bx, t) \omega(\bx,t)  \frac{\partial}{\partial t}  \int_{\R^d} G(\bx, t | \bx', 0) p_0(\bx') d\bx' d\bx dt.
\end{equation}
Shifting the time derivative of $G$ into the integral and rearranging the order of integration, \eref{eq:adjLHS} becomes
\begin{equation} 
(\diamond)=\int_{\R^d} p_0(\bx') \int_{\Gamma}  \underbrace{\int_0^T g(\bx, t) \omega(\bx,t) \frac{\partial}{\partial t} G(\bx, t | \bx', 0) dt}_{(\star)} d\bx  d\bx'.
\end{equation}
Performing integration by parts in $t$ on ($\star$), we further obtain
\numparts 
\begin{eqnarray}
\fl \label{eq:adjLHSa}  (\diamond) &= -\int_{\R^d} p_0(\bx') \int_{\Gamma}  \int_0^T G(\bx, t | \bx', 0) \frac{\partial}{\partial t}  \left(g(\bx, t) \omega(\bx,t)  \right) dt d\bx d\bx' \\
\fl \label{eq:adjLHSb} &+ \int_{\R^d} p_0(\bx') \int_{\Gamma}   G(\bx, T | \bx', 0) g(\bx, T) \omega(\bx,T) -  G(\bx, 0 | \bx', 0) g(\bx, 0) \omega(\bx,0)  d\bx d\bx' 
\end{eqnarray}
\endnumparts 
The window function $\omega \in C_0^\infty(\Gamma \times (0,T))$ models which spatio-temporal domain is accessible to the  sensor array (cf. Section \ref{subsec:MathModeling}). As the detection time is restricted, $\omega(\cdot,0) = \omega(\cdot,T) = 0$, and the second term \eref{eq:adjLHSb} vanishes. The next important step is to use the basic physical invariances of the Green's function \cite{Ba89},
\numparts
\begin{eqnarray}
\fl \hspace{4em} \label{eq:Gsym}  G(\bx, t | \bx', t') =&  G(\bx', -t' | \bx, -t), & \hspace{1.3em}\text{\small (time reversal invariance)}\\
\fl \hspace{4em} \label{eq:Gshift}  G(\bx, t | \bx', t') =&  G(\bx, t+T | \bx', t' + T), \hspace{2.5em} &\text{\small (time translation invariance)}
\end{eqnarray}
\endnumparts 
to rewrite \eref{eq:adjLHSa} as
\begin{equation}\label{eq:adjLHSa1}
\fl \hspace{3em} (\diamond) = -\int_{\R^d} p_0(\bx') \int_{\Gamma}  \int_0^T G(\bx', T | \bx, T-t) \frac{\partial}{\partial t} \left(  g(\bx, t) \omega(\bx,t) \right)  dt d\bx d\bx'.
\end{equation}  
Applying a transformation of variables $t \leftarrow T-t$, \eref{eq:adjLHSa1} becomes 
\begin{equation*}
\fl \hspace{3em} (\diamond) = \int_{\R^d} p_0(\bx') \int_{\Gamma}  \int_0^T G(\bx', T | \bx, t) \frac{\partial}{\partial t} \left( g(\bx, T- t) \omega(\bx,T -t\right)   dt d\bx d\bx'.
\end{equation*}  
If we define 
\begin{equation}
\fl \hspace{3em} \Pat^*[g](\bx') = \int_{\Gamma}  \int_0^T  G(\bx', T | \bx, t) \frac{\partial}{\partial t} \left( g(\bx, T -t) \omega(\bx,T-t) \right) dt d\bx, \label{eq:adjDef}
\end{equation} 
and compare with \eref{eq:adj}, we found the adjoint operator $\Pat^*: C_0^\infty\left(\Gamma \times (0,T)\right) \rightarrow C_0^\infty (\R^d )$ (note that $\Pat^*[g](\bx')$ is compactly supported due to the compact spatial support of $\omega$, the finite propagation time and the upper bound on $c_0(\bx)$). 
With the above derivation we can see that the adjoint can be defined through $\Pat^*[g](\bx') \mydef q(\bx',T)$, i.e. $\Pat^*$ maps $g$ to the function $q(\bx',t)$ evaluated at $t = T$, where $q$ is given by the following wave equation:
\numparts 
\begin{eqnarray} 
\fl \hspace{4em} \square^2 q(\bx',t)& =  \cases{\frac{\partial}{\partial t} \left( g(\bx', T-t) \omega(\bx',T-t) \right) \hspace{1em} \text{in} \; \Gamma \times (0, T) \\
0 \hspace{12.8em} \text{everywhere else}
}
\label{eq:apata} \\
\fl \hspace{4em} q(\bx',0)& = \hspace{1.1em} 0 \hspace{12.7em} \text{in} \quad {\mathbb{R}}^d, \label{eq:apatp0} \\
\fl \hspace{4em} \frac{\partial q}{\partial t} (\bx',0)& = \hspace{1.1em} 0 \hspace{12.7em} \text{in} \quad {\mathbb{R}}^d. \label{eq:apatpt0}
\end{eqnarray}
\endnumparts
One way to compare the adjoint operator defined by the above system (\ref{eq:apata}-\ref{eq:apatpt0}) to time reversal approaches (cf. Section \ref{subsec:ChallengesPAT}), is to choose a closed, sufficiently smooth $d-1$ dimensional surface $\Theta$ inside the sensor volume $\Gamma$, on which we can impose a boundary condition. $\Pat^{\triangleleft}$ can then be expressed as mapping $g$ to the solution $q(\bx',t)$ of the following wave equation, 
\numparts 
\begin{eqnarray}
\fl \hspace{4em} \square^2 q(\bx',t)& = 0 \label{eq:TRa} & \text{in} \quad {\mathbb{R}}^d \times (0, \infty) \setminus \Theta \times (0, T), \\
\fl \hspace{4em} q(\bx',t)& = g(\bx', T-t) \omega(\bx',T-t) \hspace{2em} & \text{on} \quad \Theta \times (0, T), \label{eq:TRb} \\
\fl \hspace{4em} q(\bx',0)& = 0 & \text{in} \quad \R^d \setminus \Theta, \label{eq:TRc} \\
\fl \hspace{4em} \frac{\partial q}{\partial t} (\bx',0)& = 0& \text{in} \quad \R^d, \label{eq:TRd}
\end{eqnarray}
\endnumparts
again, evaluated at $t = T$: $\Pat^{\triangleleft}[g](\bx') \mydef q(\bx',T)$. The main difference between time reversal and the adjoint is therefore given by how the time reversed pressure in the sensor element, $g(\bx', T-t) \omega(\bx',T-t)$, is introduced into the wave equation \eref{eq:WaveSecondOrder}: In the adjoint (\ref{eq:apata}-\ref{eq:apatpt0}), it enters as the mass source term $s(\bx,t)$ while in time reversal (\ref{eq:TRa}-\ref{eq:TRd}), it enters as an explicit constraint on a sub-set.

\section{Methods and Implementation} \label{sec:MethImp}

\subsection{Numerical Wave Propagation Model}
The adjoint model derived above could be implemented numerically using any of the standard techniques for modelling wave propagation in the time domain, so long as they allow to include a time-varying mass source term. One possible realization of such a scheme is given by the $k$-space time domain method \cite{MaSoLiTaNaWa01,CoKaArBe07,TrZhCo10}, which is already widely used for forward modelling in photoacoustics, e.g. in the \kwave/ Matlab Toolbox \cite{TrCo10}. In \ref{subsec:kWave}, we give a full technical description of how to use \kwave/ to implement discretised versions of the forward operator $\A = \M \Pat$ (\ref{eq:DefA1}-\ref{eq:DefA2}), the adjoint operator $\A^* = \Pat^* \M^*$ \eref{eq:adjDef} and the time reversal operator $\A^{\triangleleft} = \Pat^{\triangleleft} \M^*$ (\ref{eq:TRa}-\ref{eq:TRd}), which we will denote by $A \in \R^{L \times N}$, $A^* \in \R^{N \times L}$ and $A^{\triangleleft} \in \R^{N \times L}$, respectively. Note that for the discretisation, we restrict the spatial domain to a $d$-orthotope $\Omega$ (\textit{bounding box}) with $\Gamma \subset \Omega$, which is enclosed by a \termabb{Perfectly Matched Layer}{PML} to model wave propagation on an unbounded domain. 


\subsection{Image Reconstruction Methods} \label{subsec:InverseMethods}

We will now recapitulate some image reconstruction methods that rely on the discrete PAT operators described in the previous section. The first two elementary methods are obtained by directly applying the time reversal and the adjoint operator to the data:
\begin{eqnarray}
\fl \text{(TR)} \hspace{5em} p_{\text{\tiny TR}} &\mydef A^{\triangleleft} f  \label{eq:TRnum}\\
\fl \text{(BP)} \hspace{5em} p_{\text{\tiny BP}} &\mydef A^* f \label{eq:BP}
\end{eqnarray}
In \cite{StUh09,QiStUhZh11}, an enhancement of standard time reversal was proposed that takes the form of a Neumann series approach: The $k$-th iterate is obtained as
\begin{equation}
\fl \text{(iTR)} \hspace{5em} p^k_{\text{\tiny iTR}} = \sum_{j=0}^k K^j A^{\triangleleft} f, \qquad \text{with} \quad K = I_N - A^{\triangleleft} A, \label{eq:iTR}
\end{equation}
where $I_N$ is the identity in $\R^N$. A simple reformulation of \eref{eq:iTR} leads to an iterative image reconstruction algorithm:
\begin{eqnarray}
\fl  p^{k+1}_{\text{\tiny iTR}} = \sum_{j=0}^{k+1} K^j A^{\triangleleft} f &= K \sum_{j=0}^{k} K^j A^{\triangleleft} f  + A^{\triangleleft} f = K  p^{k}_{\text{\tiny iTR}} + A^{\triangleleft} f = \left( I_N - A^{\triangleleft} A\right)   p^{k}_{\text{\tiny iTR}} + A^{\triangleleft} f \nonumber\\
&\Longrightarrow \quad  p^{k+1}_{\text{\tiny iTR}} =  p^{k}_{\text{\tiny iTR}} - A^{\triangleleft} \left(A p^k_{\text{\tiny iTR}} -  f \right) \label{eq:iTRiter}
\end{eqnarray}
Note that in \cite{StUh09}, a modified time reversal operator is used that replaces \eref{eq:TRc} by a harmonic extension of the pressure field $g(\bx',T)$ to $\R^d$, which is particularly important if not all waves have left the domain, i.e., $p(\bx,T) \neq 0$ for $\bx \in \Gamma$. However, in the scenarios we consider in the numerical simulations, $T$ is chosen large enough so that the power of the residual waves is well below the noise level of the sensors. We therefore neglect the harmonic extension in this study.\\
Next, we discuss iterative image reconstruction methods based on a variational model \eref{eq:VarReg}. The simplest case is given by omitting $\mathcal{J}(p)$, i.e., we are just looking for the least-squares solution 
\begin{equation}
\fl \text{(LS)} \hspace{5em} p_{\text{\tiny LS}} \mydef \argminsub{p} \left\lbrace \case{1}{2} \sqnorm{A p -f}  \right\rbrace. \label{eq:LS}
\end{equation}
In principle, $p_{\text{\tiny LS}}$ is given as the solution of the normal equations, which we could iteratively solve with customized variants of the \emph{conjugate gradient} scheme, such as the \termabb{conjugate gradient least-squares}{CGLS} algorithm \cite{Bj96}. Equipped with a suitable stopping criterion, it yields a well known and understood \emph{regularization strategy} for solving \eref{eq:FwdEq}, i.e., a scheme that is robust to noise and ill-conditioning (cf. \cite{Ha95,CaReSh03} and the references in \cite{ArBeHa14}). However, for a better comparison to the previously presented iTR, we solve \eref{eq:LS} by a \emph{first order method} \cite{BuSaSt14}.  We only use the gradient of the cost function at $p$,
\begin{equation}
\nabla \case{1}{2} \sqnorm{A p -f}  = A^* \left(A p -f \right), 
\end{equation}
in a simple \emph{gradient descent algorithm}: For a given initialization $p^0$, the iteration is given as 
\begin{equation}
p^{k+1}_{\text{\tiny LS}} = p^k_{\text{\tiny LS}} - \eta A^* \left(A p^k_{\text{\tiny LS}} -f \right), \label{eq:GradDesc}
\end{equation}
which is very similar to \eref{eq:iTRiter}: $A^{\triangleleft}$ has been replaced by $A^*$ and we introduced a step size $\eta$.\\
As a slight modification of the above, we can incorporate the a-priori knowledge $p_0 \geqslant 0$ and solve the positivity-constrained least-squares problem
\begin{equation}
\fl \text{(LS+)} \hspace{5em} p_{\text{\tiny LS+}} \mydef \argminsub{p \geqslant 0} \left\lbrace \case{1}{2} \sqnorm{A p -f} \right\rbrace,  \label{eq:LS+}
\end{equation}
by modifying \eref{eq:GradDesc} to a \emph{gradient descent re-projection algorithm}: 
\begin{equation}
p^{k+1}_{\text{\tiny LS+}} = \Pi_+ \left(p^k_{\text{\tiny LS$_+$}} - \eta A^* \left(A p^k_{\text{\tiny LS$_+$}} -f \right) \right), \label{eq:ProGradDesc}
\end{equation}
where $\Pi_+$ is a projection onto $\R_+^N$.\\
Given the similarity between \eref{eq:iTRiter} and \eref{eq:GradDesc}, it seems natural to define a positivity-constrained version of iterative time reversal as 
\begin{equation}
 p^{k+1}_{\text{\tiny iTR+}} =  \Pi_+ \left( p^{k}_{\text{\tiny iTR$_+$}} - A^{\triangleleft} \left(A p^k_{\text{\tiny iTR+}} -  f \right) \right). \label{eq:iTR+}
\end{equation}
The next variational method employs the \termabb{total variation}{TV} energy \cite{BuOs13} as a regularizing functional:  
\begin{equation}
\fl \text{(TV+)} \hspace{5em} p_{\text{\tiny TV+}} \mydef \argminsub{p \geqslant 0} \left\lbrace \case{1}{2} \sqnorm{A p -f} + \lambda \text{TV}(p) \right\rbrace, \label{eq:TV+}
\end{equation}
where $\lambda > 0$ is the regularization parameter and $\text{TV}(p)$ is the $\ell_1$ norm of the amplitude of the gradient field of $p$ (the details of the implementation are given in the appendix). We can solve \eref{eq:TV+}, by  modifying \eref{eq:ProGradDesc} from a projected to a \emph{proximal gradient descent}:
\begin{equation}
p^{k+1}_{\text{\tiny TV+}} = \text{prox}_{\eta \lambda,\text{TV+}} \left(p^k_{\text{\tiny TV$_+$}} - \eta A^* \left(A p^k_{\text{\tiny TV+}} -f \right) \right), \label{eq:ProxGradDesc}
\end{equation}
where the \emph{proximal operator}
\begin{equation}
\text{prox}_{\alpha,\text{TV+}}(y) \mydef \argminsub{x \geqslant 0} \left\lbrace \frac{1}{2} \sqnorm{x - y} + \alpha \text{TV}(x) \right\rbrace \label{eq:TVdenoise}
\end{equation}
is given as a positivity-constrained total variation denoising problem. For our purposes, solving \eref{eq:TVdenoise} by the algorithm described in \cite{BeTe09} is sufficiently fast. \\
We chose the schemes \eref{eq:GradDesc}, \eref{eq:ProGradDesc} and  \eref{eq:ProxGradDesc} because they are easy and intuitive to explain and provide an interesting connection to iterative time-reversal. In all the above algorithms,  $\eta \in (0,2/\theta)$ with $\theta$ being the largest singular value of $A^* A$, ensures the convergence to a minimizer of the corresponding optimization problems, see \cite{BuSaSt14}. For a given experimental PAT setting, $\theta$ only has to be computed once, which can be done by a simple power iteration. 


\section{Computational Studies} \label{sec:Studies}

\subsection{Computational Scenarios}

\begin{figure}[tb]
\centering
\subfloat[][]{\includegraphics[width = 0.48\textwidth]{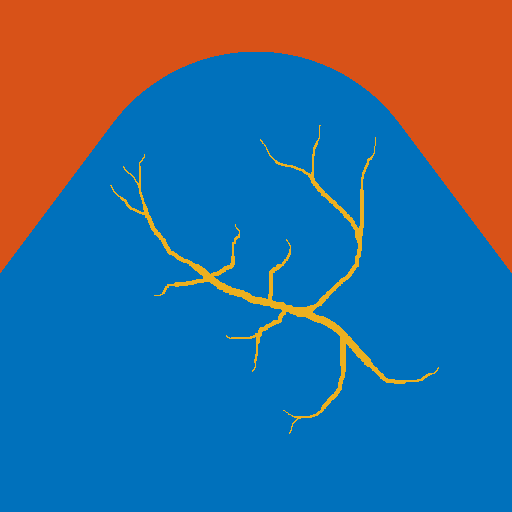} \label{subfig:2DScenarioA}}
\subfloat[][]{\includegraphics[width = 0.48\textwidth]{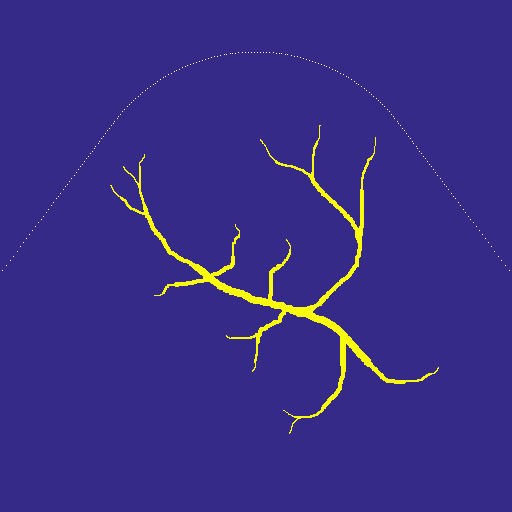} \label{subfig:2DScenarioB}}
\caption{
2D PAT scenario (II) used in the numerical studies. The spatial resolution is $N = 512^2$ plus $24$ pixels of PML layer in all directions. \protect\subref{subfig:2DScenarioA} The scenario cosists of three materials: {\textcolor{MatRed}{\textbf{Material A}}} (covering the top part of the domain): $c = 1500$, $\rho_0 = 1000$, {\textcolor{MatBlue}{\textbf{material B}}} (parabolic-like part): $c = 1400$, $\rho_0 = 1200$, {\textcolor{MatOrange}{\textbf{material C}}} (vessel-like part): $c = 1560$, $\rho_0 = 800$. \protect\subref{subfig:2DScenarioB} Ground truth $p_0$ and sensor configuration (white pixels, $200$ sensors at the interface between materials A and B.) \label{fig:2DScenario}
}
\end{figure}

To illustrate our results, we use two simple settings:
\begin{enumerate}
\item[(I)\label{sce:A}]  A 3D scenario: $\Omega = [0,1]^3$ discretised by $N = N_x^3$ isotropic voxels. The sound speed and ambient density are assumed to be homogeneous with $c = 1500 \,(m/s)$, $\rho_0 = 1000 \,(kg/m^3)$ and we add $N_{\text{\tiny PML}}$ voxels of PML at all sides of the cube. The pressure is measured in the plane $x = 0$, which corresponds to the top layer of voxels. This situation is, e.g., encountered in planar Fabry-P\'{e}rot (FP) photoacoustic scanners \cite{ZhLaPeBe09,LaJoZhTrCoPeBe12,LaNoClZhTrCoJoScLyBe12,JaLaOgTrCoZhJoPiPhMaLyPePuBe15}. 
\item[(II)\label{sce:B}]  A 2D scenario: $\Omega = [0,1]^2$ with inhomogeneous medium properties and an irregular sensor geometry, see Figure \ref{fig:2DScenario}. The medium properties vary in the range encountered, e.g., in the human breast \cite{JiWa06}.
\end{enumerate}
As our derivation in Section \ref{sec:Theory} is independent of the concrete model for $\M$, we use simple point-sampling in space and time in both scenarios.

\subsection{Precision Assessment}

\begin{figure}[t]
\centering
\includegraphics[width = \textwidth]{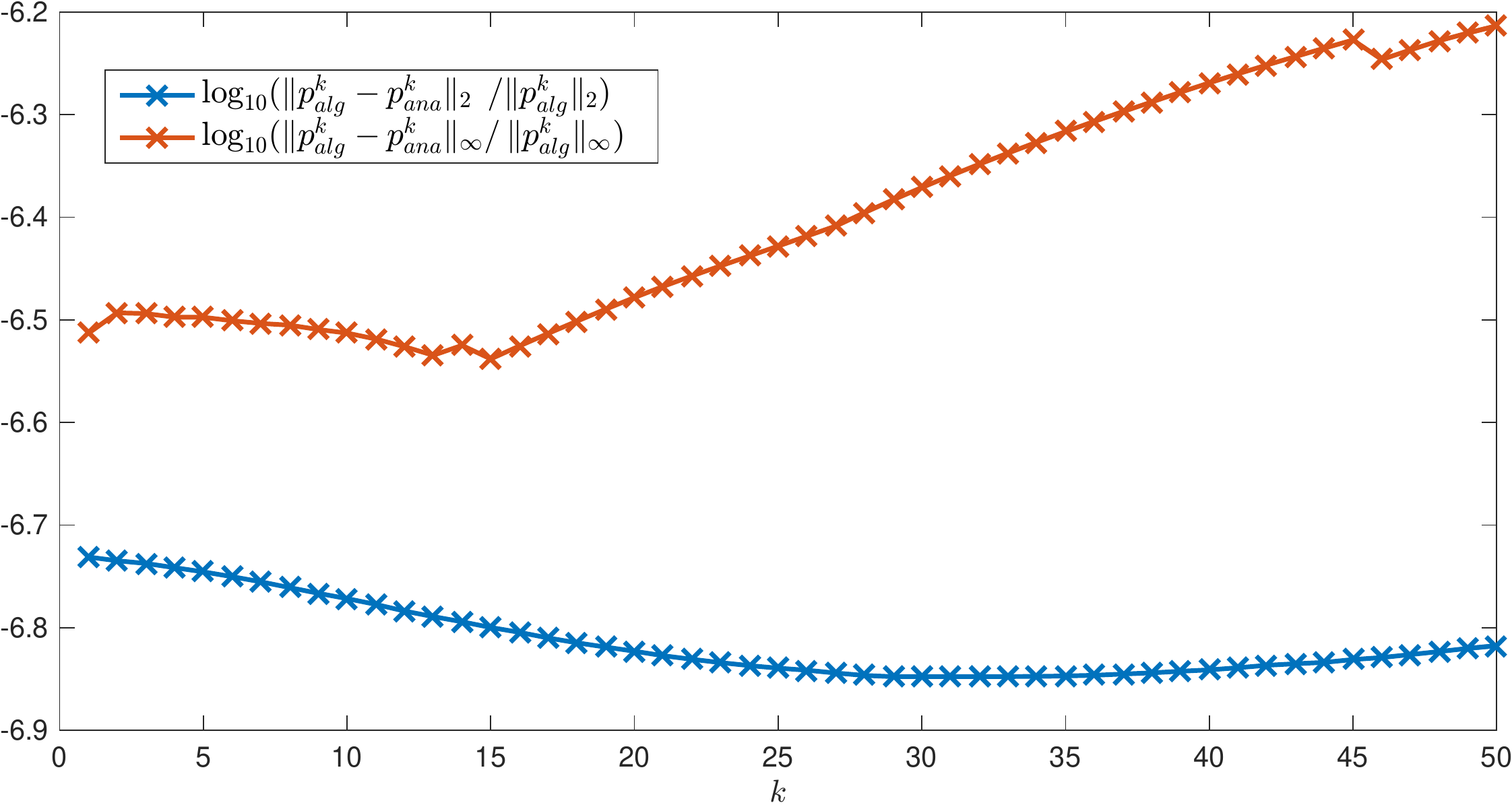}
\caption{Relative $\ell_2$ and $\ell_\infty$ errors between the iterates of scheme \eref{eq:ProGradDesc} to compute \eref{eq:LS+} using either $A^*$ or $\bar{A}^*$ as the adjoint operator. The corresponding iterates are denoted by $p_{ana}^k$ and $p_{alg}^k$, respectively.
\label{fig:Study2}}
\end{figure}

The discretised operators $A$ and $A^*$ are not fully adjoint to each other anymore, which means that if we define 
\begin{equation}
\fl \hspace{1em} \chi[A,B](x,y) \mydef | \langle A x , y \rangle - \langle x , B y \rangle |, \hspace{1em} A \in R^{L \times N}, B \in R^{N \times L}, x \in R^{N}, y \in R^{L}, \label{eq:AdjTest}
\end{equation}  
then $\chi[A,A^*](x,y) > 0$ for most $x, y$. In this section, we investigate this error by relating it to other sources of error such as the modeling error caused by the PML or numerical errors due to limited precision and examine its influence on inverse reconstructions. As a "baseline" to validate the \textit{adjoint-then-discretize} approach described in this work, we will use the \textit{discretize-then-adjoint} approach described in \cite{HuWaNiWaAn13}, which relies on the explicit adjoining of all the steps in the \kwave/ iteration (cf. \ref{subsec:kWave}). Hence, for the algebraic adjoint operator $\bar{A}^*$ obtained this way, $\chi[A,\bar A^*]$ should be 0 up to the numerical precision.\\
To test this, we computed \eref{eq:AdjTest} for $100$ random realizations of $x$ and $y$ in scenario (I). The results are summarized in Table \ref{tbl:Study1a}. To interpret the results correctly, one has to bear in mind that the PML is known to be the largest limitation of numerical accuracy of the \kwave/ iteration \cite{CoKaArBe07}, effectively preventing the use of double precision to improve numerical accuracy. Therefore, the most relevant results in Table \ref{tbl:Study1a} are the ones for $n = 128$, PML$ = 16$ and single precision. In this setting, the statistics of the error \eref{eq:AdjTest} are comparable for $A^*$ and $\bar{A}^*$. The results for double precision have to be regarded as a sanity check: While the error \eref{eq:AdjTest} for $\bar{A}^*$ is, as expected, smaller than for $A^*$, double precision is not used in practical computations. Nonetheless, we investigated the effects of this higher numerical precision in more detail in a second study using scenario (II). For this, we compare the iterates of scheme \eref{eq:ProGradDesc} to compute \eref{eq:LS+}: In the first case, we use $\bar{A}^*$ as the adjoint operator, and in the second case, we use $A^*$. The corresponding iterates will be denoted by $p_{alg}^k$ and $p_{ana}^k$, respectively. Figure \ref{fig:Study2} shows how the relative error between $p_{alg}^k$ and $p_{ana}^k$ develops with $k$. While the relative $\ell_2$ error remains almost constant, the relative $\ell_\infty$ error grows slowly but does not indicate any instability. \\
Table \ref{tbl:Study1b} completes our accuracy studies with the statistics of the $\ell_2$ and $\ell_\infty$ norm of the distance $A^*y - \bar{A}^*y$ for $100$ random realizations of $y$. We can see that in both norms the maximum/median discrepancy slightly decreases with decreasing $N$. However, we also see that this decrease can be compensated for by increasing the PML size. Furthermore, we cannot observe a difference between single and double precision in this measure anymore.\\
Summarising, the studies indicate that the proposed implementation of the analytical adjoint $A^*$ is sufficiently close to the algebraic adjoint to be confidently used in variational reconstruction schemes.

\begin{table}[t]
\begin{small}
\centering
\caption{Statistics of the error $\chi[A,B](x,y)$, $B = A^*$ (\textit{adjoint-then-discretize}) and $B = \bar{A}^*$ (\textit{discretize-then-adjoint}), for $100$ random realization of $x$ and $y$. Displayed are maximum/median of $\log_{10}$ applied to the data.}
\label{tbl:Study1a}
\textbf{single precision} \\[3pt]
\begin{tabular}{lcc}
\br
\multicolumn{3}{l}{$A^*$ (\textit{adjoint-then-discretize})} \\
 & PML $= 8$ & PML $= 16$ \\
\mr 
$N = 32^3$  & \phantom{0}-2.86/-4.55\phantom{0} & \phantom{0}-3.45/-4.99\phantom{0} \\
$N = 64^3$  & -3.12/-4.44 & -2.92/-4.60 \\
$N = 128^3$ & -2.00/-4.20 & -1.78/-4.26 \\
\br
\end{tabular} 
\hfill
\begin{tabular}{lcc}
\br
\multicolumn{3}{l}{$\bar{A}^*$ (\textit{discretize-then-adjoint})} \\
 & PML $= 8$ & PML $= 16$ \\
\mr
$N = 32^3$  & \phantom{0}-3.45/-5.14\phantom{0} & \phantom{0}-3.29/-5.14\phantom{0} \\
$N = 64^3$  & -3.30/-4.68 & -2.98/-4.63 \\
$N = 128^3$ & -2.15/-4.22 & -1.82/-4.24 \\
\br
\end{tabular} 
\vskip 10pt
\textbf{double precision} \\[3pt]
\begin{tabular}{lcc}
\br
\multicolumn{3}{l}{$A^*$ (\textit{adjoint-then-discretize})} \\
 & PML $= 8$ & PML $= 16$ \\
\mr
$N = 32^3$  &  \phantom{0}-2.79/-4.58\phantom{0} &  \phantom{0}-3.44/-5.15\phantom{0}\\
$N = 64^3$  &  -3.02/-4.72 &  -3.56/-5.28\\
$N = 128^3$ &  -2.44/-4.77 &  -2.79/-5.46\\
\br
\end{tabular} 
\hfill
\begin{tabular}{lcc}
\br
\multicolumn{3}{l}{$\bar{A}^*$ (\textit{discretize-then-adjoint})} \\
 & PML $= 8$ & PML $= 16$ \\
\mr
$N = 32^3$  &  -12.44/-13.69 & -12.12/-13.87\\
$N = 64^3$  &  -11.41/-13.47 & -11.75/-13.40\\
$N = 128^3$ &  -10.72/-12.92 & -10.97/-12.90\\
\br
\end{tabular} 
\end{small}
\end{table}

\begin{table}[t]
\begin{small}
\centering
\caption{Statistics of the $\ell_2$ and $\ell_\infty$ norm of $A^* y - \bar{A}^*y$ for $100$ random realization of $y$ (\textit{adjoint-then-discretize} vs. \textit{discretize-then-adjoint}). Displayed are maximum/median of $\log_{10}$ applied to the data.}
\label{tbl:Study1b}
\textbf{single precision} \\[3pt]
\begin{tabular}{lcc}
\br
$\ell_2$ & PML $= 8$ & PML $= 16$ \\
\mr
$N = 32^3$  & \phantom{0}-3.10/-3.13\phantom{0} &  \phantom{0}-3.79/-3.84\phantom{0}\\
$N = 64^3$  & -2.74/-2.76 &  -3.48/-3.50\\
$N = 128^3$ & -2.35/-2.36 &  -3.12/-3.13\\
\br
\end{tabular} 
\hfill
\begin{tabular}{lcc}
\br
$\ell_\infty$ & PML $= 8$ & PML $= 16$ \\
\mr
$N = 32^3$  & \phantom{0}-4.21/-4.37\phantom{0} &   \phantom{0}-4.91/-5.03\phantom{0}\\
$N = 64^3$  & -4.20/-4.32 &   -4.83/-4.97\\
$N = 128^3$ & -4.14/-4.26 &   -4.78/-4.89\\
\br
\end{tabular} 
\vskip 10pt
\textbf{double precision} \\[3pt]
\begin{tabular}{lcc}
\br
$\ell_2$ & PML $= 8$ & PML $= 16$ \\
\mr
$N = 32^3$  &  \phantom{0}-3.10/-3.13\phantom{0} &  \phantom{0}-3.79/-3.84\phantom{0}\\
$N = 64^3$  &  -2.74/-2.77 &  -3.49/-3.51\\
$N = 128^3$ &  -2.35/-2.36 &  -3.17/-3.18\\
\br
\end{tabular} 
\hfill
\begin{tabular}{lcc}
\br
$\ell_\infty$ & PML $= 8$ & PML $= 16$ \\
\mr
$N = 32^3$  &  \phantom{0}-4.21/-4.37\phantom{0} &   \phantom{0}-4.90/-5.04\phantom{0}\\
$N = 64^3$  &  -4.20/-4.32 &   -4.83/-4.97\\
$N = 128^3$ &  -4.14/-4.26 &   -4.79/-4.89\\
\br
\end{tabular} 
\end{small}
\end{table}

\subsection{Inverse Reconstructions} \label{subsec:InvRec}

\begin{figure}[tb]
\centering
\subfloat[][Ground truth $p_0$]{\includegraphics[width = 0.235\textwidth]{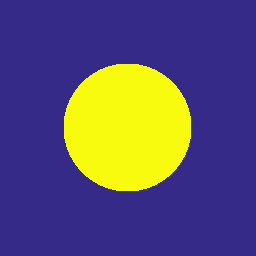} \label{subfig:3Dref}}
\subfloat[][TR, PNSR = 15.71]{\includegraphics[width = 0.235\textwidth]{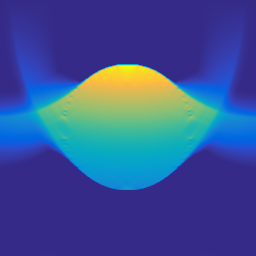} \label{subfig:3DTR}}
\subfloat[][iTR, PSNR = 18.00]{\includegraphics[width = 0.235\textwidth]{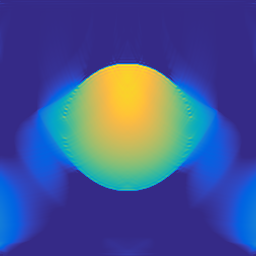} \label{subfig:3DiTR}}
\subfloat[][iTR+, PSNR = 20.97]{\includegraphics[width = 0.235\textwidth]{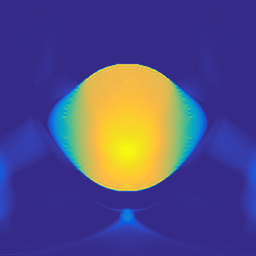} \label{subfig:3DiTR+}}\\
\subfloat[][TV+, $\lambda = 0.05$, PSNR = 23.57]{\includegraphics[width = 0.235\textwidth]{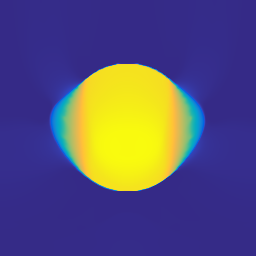} \label{subfig:3DTV}}
\subfloat[][BP, PSNR = 15.75]{\includegraphics[width = 0.235\textwidth]{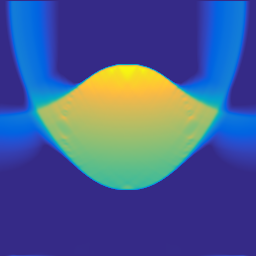} \label{subfig:3DBP}}
\subfloat[][LS, PSNR = 18.072]{\includegraphics[width = 0.235\textwidth]{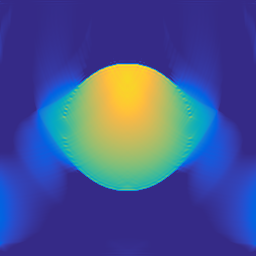} \label{subfig:3DLS}}
\subfloat[][LS+, PSNR = 21.16]{\includegraphics[width = 0.235\textwidth]{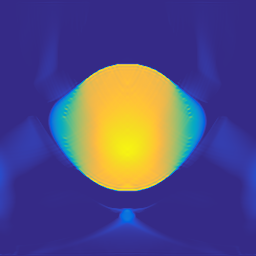} \label{subfig:3DLS+}}
\caption{2D slices at $y = 128$ through the 3D phantom $p_0$ and reconstructions thereof within scenario (I). The colorbar in each figure is scaled individually. \label{fig:3DRes}}
\end{figure}

\begin{figure}[tb]
\centering
\subfloat[][TR, PNSR = 20.99]{\includegraphics[width = 0.32\textwidth]{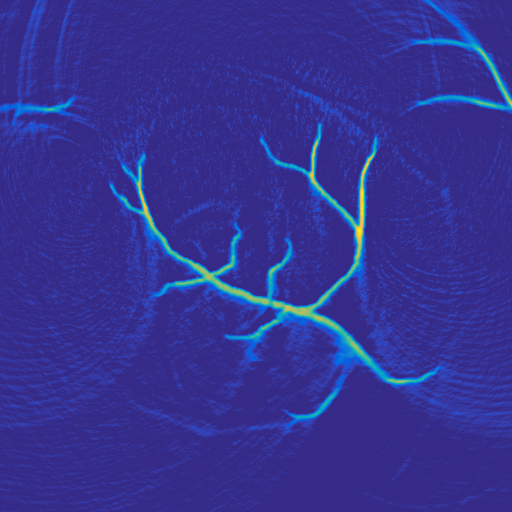} \label{subfig:2DTR}}
\subfloat[][iTR, PSNR = 23.23]{\includegraphics[width = 0.32\textwidth]{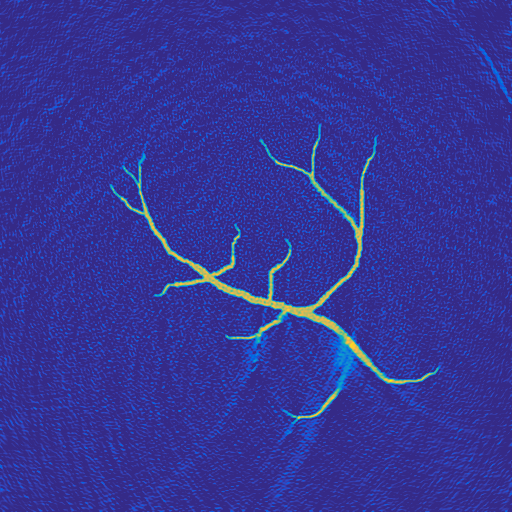} \label{subfig:2DiTR}}
\subfloat[][iTR+, PSNR = 24.61]{\includegraphics[width = 0.32\textwidth]{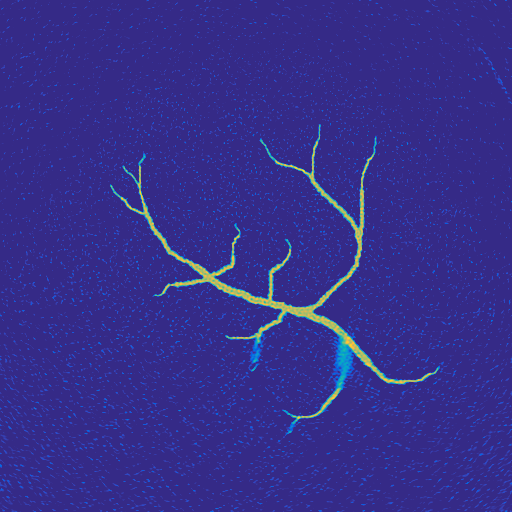} \label{subfig:2DiTR+}}\\
\subfloat[][BP, PSNR = 20.23]{\includegraphics[width = 0.32\textwidth]{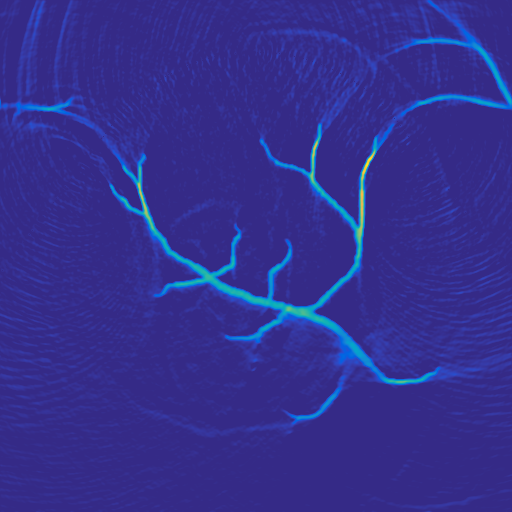} \label{subfig:2DBP}}
\subfloat[][LS, PSNR = 23.72]{\includegraphics[width = 0.32\textwidth]{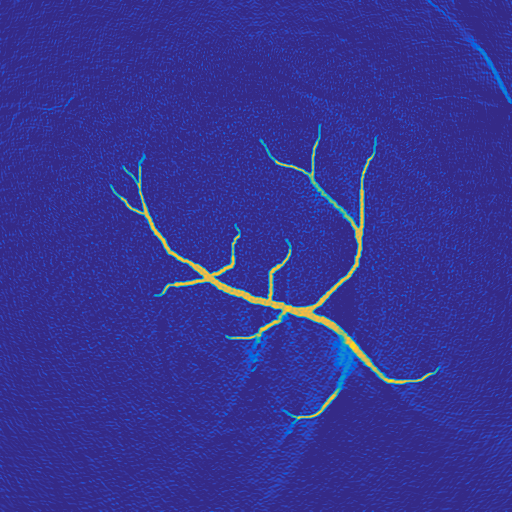} \label{subfig:2DLS}}
\subfloat[][LS+, PSNR = 24.39]{\includegraphics[width = 0.32\textwidth]{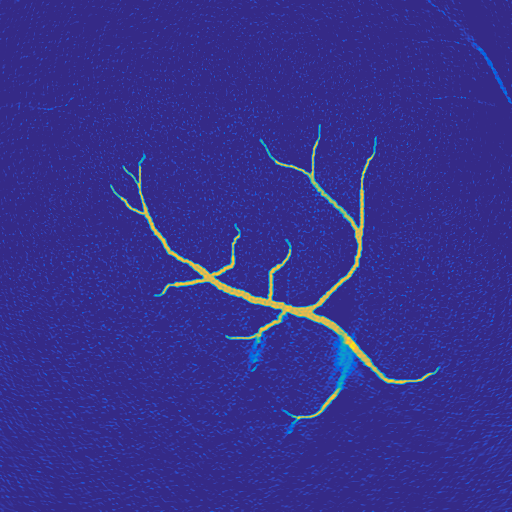} \label{subfig:2DLS+}}\\
\subfloat[][TV+, $\lambda = 0.01$, PSNR = 25.87]{\includegraphics[width = 0.32\textwidth]{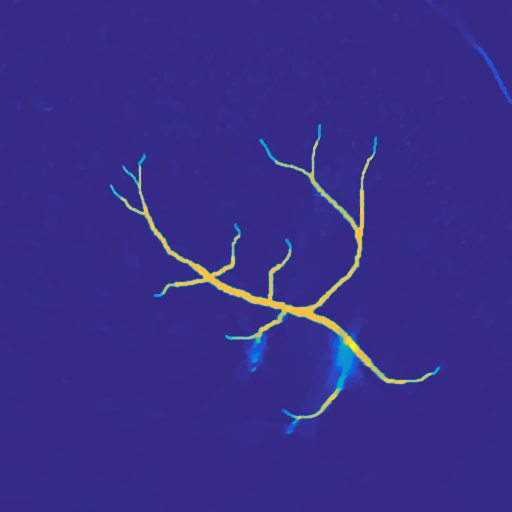} \label{subfig:2DTV}}
\caption{Results of the different image reconstruction techniques applied to scenario (II), cf. Figure \ref{fig:2DScenario}. The colorbar in each figure is scaled individually.\label{fig:2DRes}}
\end{figure}

Figure \ref{fig:3DRes} shows the results of the different inverse methods applied to scenario (I) using $N = 256^3$ voxels and a sphere as the ground truth for $p_0$. Figure \ref{fig:2DRes} shows the corresponding results for scenario (II), cf., Figure \ref{fig:2DScenario}. The \termabb{peak-signal-to-noise ratio}{PSNR} of the reconstructed solutions was computed as:
\begin{eqnarray}
\fl \text{PSNR}(p,q) = 10 \log_{10} \left(\frac{N}{\|\tilde{p}-\tilde{q}\|_2^2} \right), \qquad \text{where} \quad \tilde{p} \quad  = \text{thres}\left( \frac{p}{\|p\|_\infty},0.01\right), \\
\fl \hspace{6em} \tilde{q} \quad  \text{accordingly}, \qquad \text{and} \quad \text{thres}(x,\alpha) = \cases{ x &\text{if} \hspace{0.5em} $x \geqslant \alpha$ \\
0 &\text{else}
}
\label{eq:PSNR}
\end{eqnarray}
All computations have been performed with single precision on a GPU. For all iterative methods, $K = 100$ iterations have been used and $\eta = 1.8/\theta$ was used as a step size. All results have been projected to $\R_+^N$ as a post-processing step.\\
In line with the theoretical prediction (cf. (\ref{eq:apata}-\ref{eq:apatpt0}) and (\ref{eq:TRa}-\ref{eq:TRd})) and the numerical implementation (cf. \ref{subsec:ImplOp}), the differences between BP and TR are rather subtle. A particular reason for this is that the measurement time $T$ was chosen large enough, cf. \cite{BeGlSc15}. However, Figure \ref{fig:3DRes} shows that the limited-view geometry of scenario (I) leads to differences in the decay of intensity with depth and the appearance of the well-known circular limited-view artefacts \cite{FrQu15}. The differences between the iterative approaches based on time reversal or least-squares minimization are even more subtle than between plain TR and BP. The results rather show that including a-priori information such as positivity or total variation constraints leads to a far more substantial improvement in terms of reconstruction quality.

\section{Discussion, Conclusions and Outlook} \label{sec:DisCon}

 In this work, we examined the adjoint of the continuous PAT operator (\ref{eq:DefA1}-\ref{eq:DefA2}) from different perspectives. We first described how the adjoint is connected to the solution of a wave equation with a time-varying mass source term (\ref{eq:apata}-\ref{eq:apatpt0}) (cf. Section \ref{sec:Theory}). The simple proof was based on two main concepts: 
\begin{itemize}
\item Our mathematical model (cf. Section \ref{subsec:MathModeling}) of the measurement process explicitly accounts for the fact that the sensor array can only access the pressure field within a finite spatial volume and within a finite time interval. This was modelled by the introduction of a window function, which causes the boundary terms arising from the integration by parts in time to vanish in a central step in the proof, cf.~\eref{eq:adjLHSb}.
\item The use of the Green's function representation of the solution of the wave equation \eref{eq:WaveSecondOrder}, which reflects the temporal invariances of the underlying physical system, cf. (\ref{eq:Gsym}-\ref{eq:Gshift}). 
\end{itemize}
We would like to remark that our setting is more general in terms of choice of $\Gamma$ and $T$ than the standard one, which assumes measurements on the boundary of the domain. In this standard setting, the boundary terms arising from the integration by parts in \eref{eq:adjLHSb} only vanish under the assumption that the support of $p_0$ is contained in the interior of the domain enclosed by the measurement surface and that the time $T$ is large enough for all waves to leave this domain.\\
The explicit form of the continuous adjoint not only allows for an easy theoretical comparison to time-reversal (\ref{eq:TRa}-\ref{eq:TRd}), it also enables its numerical implementation with any existing wave propagation code that is able to include a time-varying mass source. As a particular example thereof, we examined an implementation based on the \textit{k}-space pseudospectral time domain method (cf. Sections \ref{sec:MethImp} and \ref{sec:Studies}):
\begin{itemize}
\item We demonstrated that our \textit{adjoint-then-discretize} approach is sufficiently precise compared to the \textit{discretize-then-adjoint} approach \cite{HuWaNiWaAn13} that has to be derived for each specific wave propagation code separately and yields a discrete, purely algebraic, characterization of the adjoint only.
\item We illustrated how the adjoint can be used to solve variational image reconstruction problems \eref{eq:VarReg} that incorporate a-priori constraints on the initial pressure.
\end{itemize}
An alternative description and implementation of the adjoint operator for PAT was recently presented in \cite{BeGlSc15}: In comparison to our approach, the forward problem is modelled in a more restrictive way and with a different motivation, which, by using the weak formulation approach (as opposed to the integral formulation approach we take here) leads to the characterization of the adjoint as a wave transmission problem. The latter was solved by a specific BEM-FEM approach, which is computationally intensive and the numerical studies focus on short recording times $T$ with full boundary measurements. \\
Finally, these future directions of research are directly related to this work:
\begin{itemize}
\item The inverse reconstructions (cf. Section \ref{subsec:InverseMethods} and \ref{subsec:InvRec}) revealed interesting connections and differences between iterative methods based on time reversal and the adjoint operator. This topic has to be examined in more detail, both from theoretical and numerical perspectives: As (\ref{eq:apata}-\ref{eq:apatpt0}) and (\ref{eq:TRa}-\ref{eq:TRd}) differ in the way the time-reversed pressure is introduced at the sensor elements, the effect of different sensor geometries, targets in close proximity to the sensors, artefacts caused by limited sensor coverage \cite{FrQu15} and short recording times, $T$, need to be examined. Ultimately, this is not only a mathematical question: The different artefacts caused by time reversal and adjoint operator lead to different deficiencies in in-vivo images, each causing problems for the interpretation of those images.
\item The optimization schemes \eref{eq:GradDesc}, \eref{eq:ProGradDesc} and  \eref{eq:ProxGradDesc} were only used in this work because they are easy to describe and directly relate to iterative time-reversal. However, we note here that they are computationally rather inefficient and were not tuned for optimal performance. Therefore, the number of iterations used in the numerical examples does not provide a good indicator of the expected computational cost. More sophisticated algorithms for solving the corresponding optimization problems \cite{BuSaSt14} will be examined in forthcoming work, in particular in the context of more challenging large-scale models encountered in many experimental data scenarios.
\item Code to reproduce the results of this article will be provided as example scripts for a PAT image reconstruction toolbox that will be released to complement the \kwave/ toolbox (cf. \ref{subsec:kWave}).
\end{itemize}

\ack
We gratefully acknowledge the support of NVIDIA Corporation with the donation of the Tesla K40 GPU used for this research. FL and MB acknowledge financial support from Engineering and Physical Sciences Research Council, UK (EP/K009745/1).

\appendix

\section{$k$-Space Time Domain Method} \label{subsec:kWave}

This appendix gives a technical overview of the $k$-space time domain scheme used for solution of the system of equations (\ref{eq:MonCon}-\ref{eq:PreDenRel}), as it is implemented in \kwave/\footnote{For more details, see the \kwave/ manual on {\textcolor{darkblue}{\texttt{www.k-Wave.org}}}.} \cite{TrCo10}, as well as a description of how to discretise the continuous operators defined in Section \ref{sec:Theory}. The $k$-space time domain method is a collocation scheme that interpolates between the collocation points using a truncated Fourier series. This allows field gradients to be calculated efficiently using the \termabb{fast Fourier transform}{FFT}. The `$k$-space' in the name refers to a correction $\kappa$ applied in the wavenumber domain to account for the finite difference approximation made to the time derivative. In the case of a homogeneous medium, this correction is exact, and no error arises from the temporal discretisation scheme \cite{TaMaWa02}. For acoustically heterogeneous media, the $k$-space correction still reduces the errors, e.g., those due to numerical dispersion. To solve the system (\ref{eq:MonCon}-\ref{eq:PreDenRel}) using this scheme, it is first necessary to embed the spatial domain $\Omega$ (for simplicity, we assume for the spatial dimension $d = 3$ in the following description) into a rectangular volume, which is then discretised by a regular grid of $N$ collocation points. Then, a \termabb{Perfectly Matched Layer}{PML} absorbing boundary is wrapped around the box to damp outgoing waves without reflecting them, mimicking free-space propagation. The measurement time interval $[0,T]$ is discretised by $t_n =  n \Delta t$, $n=0,\ldots,N_t$, $\Delta t = T/N_t$. The discrete pressure at iteration step $n$ is denoted by $p^n$ and the $\xi$-component of the discrete particle velocity by $u^n_\xi$. Here, $\xi$ stands for one of the spatial components $x, y, z$ of these vector fields. The terms containing $\nabla \rho_0(\bx)$ in (\ref{eq:MonCon}-\ref{eq:PreDenRel}) cancel out, and are therefore not included in the discrete scheme, to save unnecessary calculations \cite{TrJaReCo12}. The scalar `density' is split, unphysically, into three parts, $\rho^n_\xi$, to allow anisotropic absorption to be included within the PML. The ambient density will be denoted by $\tilde{\rho}_0$ and the sound speed by $c_0^2$. The \emph{k-space derivative} and the \emph{k-space operator} $\kappa$ are defined as 
\begin{equation}
\fl \hspace{3em} \frac{\partial}{\partial \xi} g \mydef \F^{-1} \left\lbrace i k_\xi \kappa \F \left\lbrace g \right\rbrace \right\rbrace; \qquad \kappa = \text{sinc}\left(\frac{c_{{\mathrm ref}} \Delta t}{2} \sqrt{k_x^2+k_y^2+k_z^2} \right), \label{eq:DisPartDev} 
\end{equation}
where $k_\xi \in \R^N$ is a the discrete wavevector in $\xi$ direction and all multiplications between $N$-dimensional vectors are understood as componentwise. $c_{\mathrm ref}$ is a homogeneous reference sound speed, chosen to ensure stability, e.g., $c_{\mathrm ref} = \max \, c_0(x)$. Staggering in time is included by interleaving the gradient and updates steps, as shown below. Grid-staggering is incorporated into the calculation of the gradients; a spatial translation of $\pm \Delta \xi$ is included as
\begin{equation}
\frac{\partial^\pm}{\partial \xi} g \mydef \F^{-1} \left\lbrace i k_\xi \kappa e^{\pm i k_\xi \Delta \xi/2} \F \left\lbrace g \right\rbrace \right\rbrace. \label{eq:DisPartDevSt} 
\end{equation}
The PML is implemented by multiplication operators $\Lambda_\xi$ and $\Lambda^s_\xi$ on the normal and staggered grid, respectively. Finally, we define $W: \R^N \longrightarrow \R^{N_\Gamma}$ to discretise the window function $\omega$ (cf. Section \ref{subsec:MathModeling}).
\paragraph{\kwave/ Iteration}
\textit{Set $p^{-1}$, $u_\xi^{-3/2}$ and $\rho_\xi^{-1}$ to zero (for all $\xi \in [x,y,z]$). Start the iteration at $n = -1$, and iterate until $N_t-1$ ($t_0 = 0$ and $t_{N_t} = T$):}
\label{eq:kWave}
\begin{eqnarray}
\fl \hspace{1em} u^{n+1/2}_\xi = \quad  \Lambda^s_\xi \left( \Lambda^s_\xi u^{n-1/2}_\xi - \frac{\Delta t}{\tilde{\rho}_0} \frac{\partial^+}{\partial \xi} p^n \right) \hspace{5.4em} \text{\small (momentum conservation)}\label{eq:IterU} \\
\fl \hspace{1em} \rho_\xi^{n+1\phantom{/2}} = \quad  \Lambda_\xi \left( \Lambda_\xi \rho_\xi^n - \Delta t \tilde{\rho}_0 \frac{\partial^-}{\partial \xi} u^{n+1/2}_\xi \right) + \Delta t \, s^{n+1/2}  \label{eq:IterRho} \hspace{2em} \text{\small (mass conservation)} \\ 
\fl \hspace{1em} p^{n+1\phantom{/2}} = \quad c_0^2 \left(\rho_x^{n+1} + \rho_y^{n+1} + \rho_z^{n+1} \right) \hspace{6.85em}  \text{\small (pressure-density relation)} \label{eq:IterP} \\
\fl \hspace{1em} g^{n+1\phantom{/2}} = \quad  W p^{n+1} \hspace{17.85em}  \text{\small (extract  pressure)} \label{eq:IterMeasure} 
\end{eqnarray}
Note that Equations \eref{eq:IterU} and \eref{eq:IterRho} are really each three equations, for $\xi \in [x,y,z]$ separately.

\section{Pseudospectral Implementation of the PAT Operators} \label{subsec:ImplOp}

Using the \kwave/ iteration, the forward operator $\Pat$ (\ref{eq:DefP1}-\ref{eq:DefP2}) and the adjoint operator  $\Pat^*$ (\ref{eq:apata}-\ref{eq:apatpt0}) can be implemented as follows:\\ 
$\Pat$ is implemented by setting 
\begin{equation} \label{eq:ImplA}
s^{n+1/2} = \frac{1}{2 c_0^2 d \Delta t } \cases{
 Q p_0&\text{for} \quad $n = -1,0$ \\
0 &\text{else}}
\end{equation}
and the result is given by $f = M g$. Here, $M$ is a discrete implementation of the measurement operator $\M$, i.e., it maps the discrete pressure time series in the sensor elements, $g \in \R^{N_t N_\Gamma}$, to the data $f \in \R^L$. In the simulation studies, $M$ is just an identity matrix. In \eref{eq:ImplA}, $Q$ is a smoothing operator introduced to eliminate the high spatial frequencies of $p_0$ to avoid unintended oscillations in the field being introduced by the band-limited interpolant (see Section 2.8 in the \kwave/ manual and \cite{TrCo11}).\\
$\Pat^*$ is implemented by setting 
\begin{equation}
\fl \hspace{3em} s^{n+1/2} = \frac{\tilde{\rho}_0}{2 d \Delta t } W^* \cases{\tilde{g}^{N_t} &for \quad $n = -1$ \\
\tilde{g}^{N_t-n+1} + \tilde{g}^{N_t-n} &for \quad $n = 0,\ldots,N_t-2$ \\
\tilde{g}^1 + 2  \tilde{g}^0 &\text{for} $\quad n = N_t-1$
}
\end{equation}
with $\tilde{g} := M^* f$ and omitting \eref{eq:IterMeasure}. The result is given by $Q p^{N_t}/(c_0^2 \tilde{\rho}_0)$, where the smoothing by the self-adjoint $Q$ has to be included to obtain adjointness with the forward implementation. The rescaling of $p^{N_t}$ and $\tilde{g}$ by $c_0^2$ and $\tilde{\rho}_0$ is necessary since we are actually solving the adjoint of the first order system (\ref{eq:MonCon}-\ref{eq:PreDenRel}) instead of the second order system \eref{eq:WaveSecondOrder}: $p$ and $\rho$ are essentially the same variables up to a scaling by $c_0^2$. In the first order adjoint system, the scaling is transferred from $\rho$ to $p$. The same holds for $\rho_0(\bx)$, which is absent from the second order system.\\
It is interesting to compare the implementation of the adjoint with the implementation of TR where the time-reversed pressure field in the sensor elements, $\tilde{g}^{{N_t}-n+1}$, is introduced as a \emph{Dirichlet source}: In the \kwave/ iteration, \eref{eq:IterRho} is replaced by 
\begin{equation}
\fl \hspace{2em} \rho_\xi^{n+1}  = \left(I_N - W W^* \right) \Lambda_\xi \left( \Lambda_\xi \rho_\xi^n - \Delta t \tilde{\rho}_0 \frac{\partial^-}{\partial \xi} \tilde{u}^{n+1/2} \right) + \frac{1}{d c_0^2 } W^*  \tilde{g}^{{N_t}-n}, \label{eq:IterRhoTR} \\ 
\end{equation}
and step \eref{eq:IterMeasure} is omitted. The result is given by $Q p^{T}$. In \eref{eq:IterRhoTR}, the multiplication with $\left(I_N - W W^* \right)$ sets the density at the sensor locations to $0$.\\


\section{Discrete Total Variation Functional} \label{subsec:TV}

If we index the pixels of a 2D image $p \in \R^N$, with $N = N_x \times N_y$ by $(i,j)$, $i = 1,\ldots,N_x$, $j = 1,\ldots,N_y$, a commonly used discretization of the total variation seminorm with Neumann boundary conditions relying on forward finite differences is given by   
\begin{equation}
\text{TV}(p) = \sum_{(i,j)}^N \sqrt{(p_{(i+1,j)}-p_{(i,j)})^2 + (p_{(i,j+1)} - p_{(i,j)})^2 }, 
\end{equation}
where $p_{(N_x+1,j)} := p_{(N_x,j)}$, $p_{(i,N_y+1)} := p_{(i,N_y)}$. The implementation in 3D works in the same way.


\noappendix 

\section*{References}

\bibliographystyle{abbrv}
\bibliography{AdjointBibfile}

\end{document}